# Tight Chromatic Upper Bound for {$3K_1$, $2K_1+(K_2\cup K_1)$}-free Graphs

Medha Dhurandhar
mdhurandhar@gmail.com

**Abstract**

Problem of finding an optimal upper bound for the chromatic no. of $3K_1$-free graphs is still open and pretty hard. It was proved by Choudum et al that upper bound on the chromatic no. of {$3K_1$, $2K_1+(K_2\cup K_1)$}-free graphs is $2\omega$. We improve this by proving for such graphs $\chi \leq 8 = \lceil \frac{3\omega}{2} \rceil$ for $\omega = 5$ and $\chi \leq \frac{3\omega}{2}$ for $\omega \neq 5$, where $\omega$ is the size of a maximum clique in G. We also give examples of extremal graphs.

**1. Introduction:**

In [1], [2], [4], [5] chromatic bounds for graphs are considered especially in relation with $\omega$ and $\Delta$. Gyárfás [6] and Kim [7] show that the optimal $\chi$-binding function for the class of graphs with independence number 2 has order $\omega^2/\log(\omega)$. Independence number 2 implies that the graph is $3K_1$-free. If we forbid additional induced subgraphs, the order of the optimal $\chi$-binding function may drop below $\omega^2/\log(\omega)$. We consider in this paper only finite, simple, connected, undirected graphs. The vertex set of G is denoted by V(G), the edge set by E(G), the maximum degree of vertices in G by $\Delta$(G), the maximum clique size by $\omega$(G) and the chromatic number by $\chi$(G). N(u) denotes the neighbourhood of u and $\overline{N(u)}$ = N(u)+u. In [1] it was proved that is $2\omega$ is an upper bound on the chromatic no. of {$3K_1$, $2K_1+(K_2\cup K_1)$}-free graphs and the problem of finding tight chromatic upper bound for a {$3K_1$, $2K_1+(K_2\cup K_1)$}-free graph was stated as open.

For further notation please refer to Harary [3].

**2. Main Result:**

Before proving the main result we prove some Lemmas.

**Lemma 1:** Let G be {$3K_1$, $2K_1+(K_2\cup K_1)$}-free. Then $V(G) = \bigcup_{1}^{j} M_i$ where $1\leq j\leq 4$ and

**1.1** $<M_i>$ is complete $\forall$ i, $1\leq i\leq j$.
**1.2** Every vertex of $M_2$ is non-adjacent to exactly one vertex of $M_1$.
**1.3** If $m_i \in M_i$ for i = 1, 2 are s.t. $m_1m_2 \notin E(G)$ then every $m_j \in M_j$ (j=3, 4) is adjacent to exactly one of $m_1$ and $m_2$.
**1.4** Every x, y $\in M_i$ (i=3, 4) have at least $|M_2|$-2 common adjacent vertices from $M_1 \cup M_2$.
**1.5** If b $\in M_3$, c $\in M_4$ are s.t. bc $\in$ E(G), then b and c have at least $|M_2|$-1 common adjacent vertices from $M_1 \cup M_2$.
**1.6** If $|M_1| \geq |M_2| \geq 4$, then either $<M_3\cup M_4>$ is complete or $m_im_j \notin E(G)$ for $3\leq i, j\leq 4$.
**1.7** Let b $\in M_3$ and c, c' $\in M_4$ be s.t. bc, bc' $\in$ E(G). If m $\in M_1\cup M_2$ is s.t. cm, c'm $\in$ E(G) (cm, c'm $\notin$ E(G)), then bm $\in$ E(G) (bm $\notin$ E(G)).

Proof: If G is complete, then the **Lemma** is trivially true. Let $\exists$ v, w $\in$ V(G) s.t. vw $\notin$ E(G). Let A = {x $\in$ V(G)/xv, xw $\in$ E(G)}, B = {x $\in$ V(G)/xv $\in$ E(G), xw $\notin$ E(G)}, and C = {x $\in$ V(G)/xw $\in$ E(G), xv $\notin$ E(G)}. Further let D $\subseteq$ A be s.t. <D> is a maximum clique in <A> and Y = A–D.

**1.1** Let D = $M_1$, Y = $M_2$, B = $M_3$ and C = $M_4$. As G is $3K_1$-free, <B>, <C> are complete. Let y, z $\in$ Y. Let if possible yz $\notin$ E(G). Let x $\in$ D be s.t yx $\notin$ E(G). As G is $3K_1$-free zx $\in$ E(G) and <v, z, w, x, y> = $2K_1+(K_2 \cup K_1)$, a contradiction.

**1.2** By definition every vertex of Y is non-adjacent to at least one vertex of D. Let if possible y ∈ Y be non-adjacent to x, z ∈ D. Then <v, x, w, z, y> = $2K_1+(K_2 \cup K_1)$, a contradiction.

**1.3** Let x ∈ D and y ∈ Y be s.t. xy ∉ E(G). As G is $3K_1$-free, every b ∈ B (c ∈ C) is adjacent to x or y. Let if possible bx, by ∈ E(G). But then <x, b, y, v, w> = $2K_1+(K_2 \cup K_1)$, a contradiction.

**1.4** Let Y' ⊆ D be s.t. a vertex of Y is non-adjacent to some vertex of Y'. Let if possible b, b' ∈ B have at most |Y|-3 common adjacent vertices from Y∪Y'. W.l.g. let $y_i$' ∈ Y' and $y_i$ ∈ Y be s.t. $y_i y_i$' ∉ E(G) for 1≤i≤3, and $by_i$, $b'y_i$' ∈ E(G). But then <b, $y_1$, $y_3$', $y_2$, b'> = $2K_1+(K_2 \cup K_1)$, a contradiction.

Let X = D-Y', Y = $\bigcup_i y_i$, Y' = $\bigcup_i y_i'$ with $y_i' y_i$ ∉ E(G) ∀ i.

**1.5** Next let b ∈ B, c ∈ C be s.t. bc ∈ E(G). Let if possible b, c have at the most |Y|-2 common adjacent vertices from Y∪Y'. W.l.g. let $by_i$', $cy_i$ ∈ E(G) for 1≤i≤2. But then <b, $y_1$', w, $y_2$', c> = $2K_1+(K_2 \cup K_1)$, a contradiction.

**1.6** W.l.g. let if possible ∃ b ∈ B; c, c' ∈ C s.t. bc ∈ E(G) and bc' ∉ E(G) and $by_i$ ∈ E(G) ∀ i. As G is $3K_1$-free, $c'y_i$' ∈ E(G) ∀ i. As |Y| ≥ 4, by **1.5**, $cy_i$ ∈ E(G) for at least three i, a contradiction to **1.4**.

**1.7** If bm ∉ E(G), then <m, c, b, c', v> = $2K_1+(K_2 \cup K_1)$, a contradiction.

This proves **Lemma 1**.

**Note: Lemma 1** defines the structure of {$3K_1$, $2K_1+(K_2 \cup K_1)$}-free graphs.

**Lemma 2**: If G is $3K_1$-free and ω ≤ 3, then $\chi(G) \leq \frac{3\omega}{2}$.

Proof: The result is trivially true for ω ≤ 2. Let ω = 3. Then Δ ≤ 5 and |V(G)| ≤ 8. If Δ ≤ 4, then $\chi(G) \leq \lceil \frac{\Delta+\omega+1}{2} \rceil \leq 4 < \frac{3\omega}{2}$. Hence let Δ = 5, V(G) = $\bigcup_i a_i$ and deg $a_1$ = 5. Then N($a_1$) = <$a_2$, $a_3$, $a_4$, $a_5$, $a_6$> = $C_5$. If |V(G)| < 8, then clearly $\chi(G)$ = 4. Hence let |V(G)| = 8. As ω=3, $a_i a_j$ ∉ E(G) for some i ∈ {2, 3}, j ∈ {7, 8}. W.l.g. let $a_2 a_7$ ∉ E(G). Then color pairs ($a_1$, $a_8$), ($a_2$, $a_7$), ($a_3$, $a_5$), ($a_4$, $a_6$) by same colors to get a 4-coloring of G.

This proves **Lemma 2**.

**Definition:** We call the pair (v, w) described in **Lemma** 1, a **partitioning pair**. Henceforth let X = $\bigcup_i x_i$; Y = $\bigcup_i y_i$, Y' = $\bigcup_i y_i'$ with $y_i' y_i$ ∉ E(G); B = $\bigcup_i b_i$ and C = $\bigcup_i c_i$ be as defined in **Lemma 1.**

**Lemma 3:** If G is {$3K_1$, $2K_1+(K_2 \cup K_1)$}-free, then **(A)** ω = 4 ⟹ $\chi(G) \leq 6$ **(B)** ω = 5 ⟹ $\chi(G) \leq 8$.
Proof: Let if possible G be a smallest {$3K_1$, $2K_1+(K_2 \cup K_1)$}-free graph not satisfying the result. Let deg v = Δ, (v, w) be a **partitioning pair**. Then every u ∈ V(G) is non-adjacent to at least |C|+1 but not more than ω vertices in G. Also as <X∪Y∪v>, <B∪v>, <C∪w> are complete, |X∪Y| ≤ ω-1, |C| ≤ |B| ≤ ω-1. Also |X∪Y∪Y'| ≥ Δ-ω+1                **I**

**(A)**   Let ω = 4. Then Δ ≤ 8 and If Δ ≤ 7, then $\chi \leq \lceil \frac{\Delta+\omega+1}{2} \rceil$ = 6. Hence let Δ = 8. Then by **I**, |X∪Y∪Y'| ≥ 5 ⟹ |X| ≤ 3 and 2 ≤ |Y| = |Y'| ≤ 3.

First let |Y| = 3. Then as every $b_i$ ($c_j$) is non-adjacent to three vertices from Y∪Y' and w (v),

$<B \cup C>$ is complete. Also as $\omega = 4$, by **I**, $|X| = 0$ and as $\Delta = 8$, $|B| = 2 \Rightarrow |C| \leq 2$. If $|C| \leq 1$, then color the pairs $(y_i', y_i)$ $(1 \leq i \leq 3)$; $(v, c_1)$; $(w, b_1)$; and $b_2$ by a new color to get $\chi = 6$, a contradiction. Hence let $|C| = 2 \Rightarrow y_i'$ $(y_i)$ is non-adjacent to at least two but not more than three vertices of $B \cup C$. Thus by **1.7** of **Lemma 1**, every $y_i'$ $(y_i)$ is non-adjacent to one $b_i$ and one $c_k \Rightarrow$ W.l.g. let $c_1 y_i' \in E(G)$ $\forall$ i. Then by **1.5** of **Lemma 1**, w.l.g. $y_1' b_i \in E(G)$ $\forall$ i = 1, 2, a contradiction.

Next let $|Y| = 2$. Then $|X| = 1$, $|B \cup C| = 3$ and hence $|B| \geq 2$. If $|B| = 2$, then color the pairs $(y_i, y_i')$ $(1 \leq i \leq 2)$; $(w, b_1)$; $(v, c_1)$ by same colors and x, $b_2$ by two new colors to get $\chi \leq 6$, a contradiction. Hence $|B| = 3$. W.l.g. let $xb_1 \notin E(G)$. Then color the pairs $(x, b_1)$; $(v, w)$; $(y_i, y_i')$ $(1 \leq i \leq 2)$ by same colors and $b_2, b_3$ by two new colors to get $\chi \leq 6$, a contradiction.

**(B)** Let $\omega = 5$. Then $\Delta \leq 13$ and $9 \leq \chi \leq \lceil \frac{\Delta+6}{2} \rceil \Rightarrow 11 \leq \Delta \leq 13$. Also $|C| \leq |B| \leq 4$, $|X \cup Y \cup Y'| \geq 7$, $|X \cup Y| \leq 4$; $|Y| \geq 3$ and $|X| \leq 1$.

If $|X| = 1$, then $|Y| = 3$ and $|B| \geq 4 \Rightarrow |B| = 4 \Rightarrow \exists$ $b_i$ s.t. $xb_i \notin E(G)$ (else $<x, v, b_1, b_2, b_3, b_4> = K_6$). W.l.g. let $xb_1 \notin E(G)$. If $|C| \leq 1$, then color pairs $(x, b_1)$; $(w, b_2)$; $(v, c_1)$; $(y_i, y_i')$ $(1 \leq i \leq 3)$ by same colors and $b_3, b_4$ by new colors to get $\chi \leq 8$, a contradiction. Hence $|C| > 1 \Rightarrow$ x is non-adjacent to at least three vertices of G. Also by **1.7** of **Lemma 1**, x is non-adjacent to at most one vertex of B and one of C, a contradiction.

Hence $|X| = 0 \Rightarrow |Y| = 4$ and $|B| \geq 3$. By **1.3** of **Lemma 1**, every $b_i$ is adjacent to 4 vertices from $Y \cup Y'$. W.l.g. let $b_1 y_i \in E(G)$ for $1 \leq i \leq 4$. But then $<b_1 \cup Y \cup v> = K_6$, a contradiction.

This proves **Lemma 3**.

**Main Result:** If G is $\{3K_1, 2K_1+(K_2 \cup K_1)\}$-free and $\omega \geq 6$, then $\chi(G) \leq \frac{3\omega}{2}$.

Proof: Let if possible G be a smallest such graph with $\chi(G) > \frac{3\omega}{2}$. Let $v \in V(G)$ be with deg $v = \Delta$ and $s = |V(G)|-\Delta-1$. Now $s \geq 1$ (else $\chi(G-v) \leq \frac{3\omega(G-v)}{2}$ and $\chi(G) \leq \chi(G-v)+1 \leq \frac{3(\omega(G)-1)}{2}+1 \leq \frac{3\omega}{2}$). Again $s \geq 2$ (else let $vw \notin E(G) \Rightarrow$ in any $(\frac{3\omega}{2})$-coloring of G-w, v receives a unique color say $\alpha$. By coloring w by $\alpha$, we get $\chi(G) \leq \frac{3\omega}{2}$). Hence $\Delta \leq |(V(G)|-3$ and as $\frac{3\omega}{2} < \chi(G) \leq \lceil \frac{\Delta+\omega+1}{2} \rceil$, $\Delta \geq 2\omega-1 \Rightarrow \Delta+s+1 \geq 2\omega+2$. Also $\forall$ $\omega \geq 6$, $\omega+3 \leq \frac{3\omega}{2} < \chi(G)$. **I**

Let $v, w \in V(G)$ be s.t. deg $v = \Delta$ and $vw \notin E(G)$. Let $(v, w)$ be a **partitioning** pair. Then $|X|+|Y|+1$, $|B|+1$, $|C|+1 \leq \omega$ and as $s \geq 2$, $1 \leq |B| \leq |C|$. Also $|X \cup Y \cup Y'| \geq \Delta-\omega+1$ and hence by **I**, $|X \cup Y \cup Y'| \geq \omega \Rightarrow |Y| > 0$. **II**

**Case 1:** $b_i c_j \notin E(G)$ $\forall$ i and j.
W.l.g. let $b_1 y_i \in E(G)$ $\forall$ i $\Rightarrow c_j y_k' \in E(G)$ $\forall$ j, k $\Rightarrow b_j y_k \in E(G)$ $\forall$ j, k. Also let $X_1 = \{x \in X/ xb_j \in E(G)$ $\forall$ j$\}$ and $X_2 = X-X_1$. Clearly $<X_1 \cup B \cup Y \cup v>$ and $<X_2 \cup C \cup Y' \cup w>$ are complete and $\Delta+s+1 = ||X_1 \cup B \cup Y \cup v|+|X_2 \cup C \cup Y' \cup w| \leq 2\omega$, a contradiction to **I**.

**Case 2:** $\exists$ $b_i, b_j, c_k, c_m$ s.t. $b_i c_k \in E(G)$ and $b_j c_m \notin E(G)$.
By **1.6** of **Lemma 1**, $|Y| \leq 3$. Let $B_1 = \{b_{1i}/b_{1i} c_j \notin E(G)$ for some j$\}$, $B_2 = B-B_1$, $C_1 = \{c_{1i}/ c_{1i} b_j \notin$

E(G) for some j} and $C_2 = C-C_1$. Clearly $|B_1|, |C_1| \geq 1$. Let $B_j = \bigcup_i b_{ji}$, $C_j = \bigcup_i c_{ji}$, j = 1, 2.

**Claim 1:** If $b_{mi}c_{nj}, b_{mi}y_1, c_{nj}y_1' \in E(G)$, then $x_pb_{mi} \in E(G)$ iff $x_pc_{nj} \in E(G)$.
If $x_pb_{mi} \in E(G)$ and $x_pc_{nj} \notin E(G)$, then $<x_p, w, c_{nj}, y_1', b_{mi}> = 2K_1+(K_2 \cup K_1)$, a contradiction.

**Claim 2:** If some $b_{1i}$ ($c_{1i}$) is non-adjacent to two vertices of $C_1$ ($B_1$) and $b_{1i}y_j \in E(G)$ $\forall$ j, then $b_{1k}y_j$, $c_{1m}y_n' \in E(G)$ $\forall$ k, j, m, n.
Let $b_{1i}c_{1k}, b_{1i}c_{1m} \notin E(G)$. Then $c_{1k}y_p', c_{1m}y_p' \in E(G)$ $\forall$ p. If $\exists$ $b_{1n}$ s.t. say $y_1b_{1n} \notin E(G)$, then $y_1'b_{1n} \in E(G)$ and $b_{1n}c_{1k}, b_{1n}c_{1m} \in E(G) \Rightarrow \exists$ $c_{1t}$ s.t. $c_{1t}b_{1n} \notin E(G) \Rightarrow c_{1t}y_1 \in E(G)$ and $c_{1t}b_{1i} \in E(G)$. But then $<b_{1n}, c_{1k}, c_{1t}, c_{1m}, b_{1i}> = 2K_1+(K_2 \cup K_1)$, a contradiction.

**Claim 3:** If every $b_{1i}$ is non-adjacent to exactly one $c_{1i}$ and vice versa, then $\forall$ j, $x_j$ is adjacent to exactly one of $b_{1i}$ and $c_{1i}$ where $b_{1i}c_{1i} \notin E(G)$.
Let if possible $x_jb_{1i}, x_jc_{1i} \in E(G)$. Let $b \in B-b_{1i}$. Then as $c_{1i}b \in E(G)$, $x_jb \in E(G)$ (else $<x_j, b_{1i}, b, v, c_{1i}> = 2K_1+(K_2\cup K_1)$). Similarly $\forall$ $c \in C-c_{1i}$, $x_jc \in E(G)$, a contradiction to $s \geq 2$.

**Claim 4:** If $|Y| \geq 2$, then every $b_{1i}$ is non-adjacent to exactly one $c_{1j}$ and vice versa.
If not, then $|B_1| \geq 2$ or $|C_1| \geq 2$ and by **Claim 2** w.l.g. let $b_{1i}y_j, c_{1i}y_j' \in E(G)$ $\forall$ i, j. As $|Y| \geq 2$, by **1.5** of **Lemma 1**, $b_{1i}c_{1j} \notin E(G)$ $\forall$ i, j. Also $|B_1| \geq 2$ (else $|C_1| \geq 2$ and as $|B| \geq |C|$, $|B_2| > 0$ and $b_{2i}y_j' \in E(G)$ $\forall$ i, j $\Rightarrow <b_{21}, b_{11}, y_1, b_{12}, y_2'> = 2K_1+(K_2\cup K_1)$). Let $X_1 = \{x_i / x_ib_{1j} \notin E(G)$ for some j} and $X_2 = X-X_1$. Then $|B_2 \cup C_2| \geq 2$ (else $\Delta+s+1 = |X_1 \cup C_1 \cup Y' \cup w|+|X_2 \cup B_1 \cup Y \cup v|+|B_2 \cup C_2| \leq 2\omega+1$). Also $|B_2| \geq 1$ (else $|C_2| \geq 2$, by **1.7** of **Lemma 1**, $c_{2i}y_j \in E(G)$ $\forall$ i, j and $<c_{11}, c_{21}, y_1, c_{22}, y_2'> = 2K_1+(K_2\cup K_1)$). Now $|Y| = 2$ (else by **1.5** of **Lemma 1**, $b_{21}$ is adjacent to at least two vertices of Y' say $y_1', y_2'$ and $<b_{21}, b_{11}, y_1, b_{12}, y_2'> = 2K_1+(K_2\cup K_1)$). Then by **II**, $|X| \geq 2$. Also $b_{2i}$ is adjacent to one vertex of Y' by **1.5** of **Lemma 1** and one vertex of Y (else $<b_{21}, b_{11}, y_1, b_{12}, y_2'> = 2K_1+(K_2\cup K_1)$). W.l.g. let $y_1'b_{21}, y_2b_{21} \in E(G)$. Now $x_ic_{1j} \in E(G)$ $\forall$ i, j (else $x_ib_{1j} \in E(G)$ $\forall$ j $\Rightarrow$ By **Claim 1**, $x_ib_{21} \notin E(G) \Rightarrow <x_i, b_{11}, b_{21}, b_{12}, y_1'> = 2K_1+(K_2\cup K_1)$). Again by **Claim 1**, $x_ib_{2i} \in E(G)$ $\forall$ i and if $|C_2| > 0$, then $c_{2i}y_j \in E(G)$ $\forall$ i, j and $x_ic_{2j} \in E(G)$ $\forall$ j. As $s \geq 2$, w.l.g. let $x_ib_{11}, x_ib_{12} \notin E(G)$. Then $<b_{21}, b_{11}, y_1, b_{12}, x_i> = 2K_1+(K_2\cup K_1)$, a contradiction.

**Claim 5:** $|B_2| \leq 1$ and $|C_2| \leq 1$.
Let if possible $|B_2| \geq 2$. W.l.g. let $b_{11}c_{11} \notin E(G)$ and $c_{11}y_j \in E(G)$ $\forall$ j $\Rightarrow$ $b_{11}y_j' \in E(G)$ $\forall$ j.

**A:** Let $b_{2i}y_j \in E(G)$ $\forall$ i, j $\Rightarrow$ By **1.7** of **Lemma** 1, $c_{mn}y_j \in E(G)$ $\forall$ m, n, j $\Rightarrow$ $b_{1i}y_j' \in E(G)$ $\forall$ i, j and every $b_{1i}$ is adjacent to at the most one $c_{1j}$ and vice versa, $|Y| = 1$ (else $<b_{11}, b_{21}, y_1, b_{22}, y_2'> = 2K_1+(K_2\cup K_1)$) and by **1.7** of **Lemma** 1, $|C_2| \leq 1$. Also by **II**, $|X| > 0$. As $s \geq 2$, for each $b_{2i}$ $\exists$ say $x_i$ s.t. $x_ib_{2i} \notin E(G)$. Again $x_i \neq x_j$ for $i \neq j$ (else $x_ib_{2i}, x_ib_{2j} \notin E(G)$, and if $x_ib_{11} \in E(G)$ then $<b_{11}, b_{21}, y_1, b_{22}, x_i> = 2K_1+(K_2\cup K_1)$ and if $x_ic_{11} \in E(G)$ then $<c_{11}, b_{2i}, v, b_{2j}, x_i> = 2K_1+(K_2\cup K_1)$). Again every $x_k$ is non-adjacent to some $b_{2i}$ (else as $|B_2| \geq 2$, clearly $x_kb_{1m}, x_kc_{1m} \in E(G)$ and deg $x_k > \Delta$). Then $|B_1| \leq 2$ (else if say $x_1b_{11}, x_1b_{12} \in E(G)$ then $<x_1, b_{11}, b_{21}, b_{12}, y_1> = 2K_1+(K_2\cup K_1)$ and if $x_1b_{11}, x_1b_{12} \notin E(G)$ then $<b_{22}, b_{11}, y_1', b_{12}, x_1> = 2K_1+(K_2\cup K_1)$). Similarly $|C_1| \leq 2$. W.l.g. let $b_{1i}c_{1i} \notin E(G)$ for i $\leq 2$. Then color pairs $(x_i, b_{2i})$; $(y_1, y_1')$; $(v, w)$; $(b_{1i}, c_{1i})$, and vertex of $C_2$ (if any) by a new color to get $\chi \leq |X|+|Y|+1+3 \leq \omega+3$, a contradiction.

**B:** $\exists$ $b_{21}, b_{22}$ s.t. $y_1'b_{21}, y_1b_{22} \in E(G)$. Then by **1.7** of **Lemma 1**, $|C| \leq 2$ and $y_1b_{2i} \in E(G)$ $\forall$ i > 1. Also by **1.5** of **Lemma 1**, $<b_{21} \cup \{Y-y_1\}>$ is complete.

**B.1:** $|C| = 2$ and $c \in C-c_{11} \Rightarrow y_1'c \in E(G)$ and by **1.7** of **Lemma 1**, $|B_2| = 2$. Now $|B_1| \geq 3$ (else color vertices of X by $|X|$ new colors and pairs $(y_1, b_{21})$; $(y_1', b_{22})$; those of corresponding non-adjacent vertices of Y-y and Y'-z; $(b_{11}, c_{11})$; $(w, b_{12})$; $(v, c)$ by same colors to get $\chi \leq |X|+|Y|+1+3 \leq \omega+3$).

Clearly some $c_{1i}$ is non-adjacent to more than one vertex of $B_1$ and by **Claim 2**, $|Y| = 1$ and $\langle B_1 \cup Y' \rangle$ is complete. Then by **1.7** of **Lemma 1**, $c_{11}$ is adjacent to at the most one $b_{1i}$. Also as $s \geq 2$, $\exists\ x \in X$ s.t. $xb_{22} \notin E(G) \Rightarrow x$ is adjacent to at the most one $b_{1i}$ (else $\langle x, b_{1i}, b_{22}, b_{1k}, y_1 \rangle = 2K_1+(K_2 \cup K_1)$). Also $x$ is adjacent to at the most one vertex of $C$ (else $\langle x, c_{11}, b_{22}, c, v \rangle = 2K_1+(K_2 \cup K_1)$). Clearly $xc_{11} \in E(G) \Rightarrow xc \notin E(G)$ and $\langle c \cup B_1 \rangle$ is complete $\Rightarrow c \in C_2$. Now by **II**, $|X| \geq 4$. Now $\langle X \cup c_{11} \rangle$ is complete (else let $x' \in X$ be s.t. $x'c_{11} \notin E(G) \Rightarrow x'b_{1i} \in E(G)\ \forall\ i \Rightarrow$ By **1.7** of **Lemma 1**, $x'c \in E(G) \Rightarrow \langle x, x', c, y_1', c_{11} \rangle = 2K_1+(K_2 \cup K_1)) \Rightarrow cx' \notin E(G)\ \forall\ x' \in X$ and $\langle c_{11}, x, y_1', x', c \rangle = 2K_1+(K_2 \cup K_1)$, a contradiction.

**B.2:** $|C| = |C_1| = 1$. Then $\langle Y' \cup B_1 \rangle$ is complete. If $|Y| \geq 2$, then by **Claim 4**, $|B_1| = 1$. Also $|B_2| \leq 3$ (else by **1.7** of **Lemma 1**, $mb_{2j}, mb_{2k} \in E(G)$ for some $j, k > 1$ and $\langle b_{11}, b_{2j}, y_1, b_{2k}, n \rangle = 2K_1+(K_2 \cup K_1)$). Then color vertices of $X$ and $b_{23}$ by $|X|+1$ new colors and pairs $(y_1, b_{21}); (y_1', b_{22});$ $(y_i, y_i'); (v, w)$ and $(b_{11}, c_{11})$ by same colors to get $\chi \leq |X|+|Y|+1+3| \leq \omega+3$, a contradiction. Hence $|Y| = 1$ and by **II**, $|X| \geq 4$. Also as $s > 1$, $\exists\ x \in X$ s.t. $xb_{2i} \notin E(G)\ \forall\ i$. Now $\forall\ x \in X$, either $xc_{11} \in E(G)$ or $xb_{1i} \in E(G)$, but not both (else $\exists\ x$ s.t. $xc_{11}, xb_{1i} \in E(G) \Rightarrow xb_{2k} \in E(G)\ \forall\ k$ as otherwise $\langle x, b_{1i}, b_{2k}, v, c_{11} \rangle = 2K_1+(K_2 \cup K_1) \Rightarrow$ As $s > 1$, $xb_{1t}\ xb_{1m} \notin E(G)$ and $\langle b_{22}, b_{1t}, y_1', b_{1m}, x \rangle = 2K_1+(K_2 \cup K_1)$). W.l.g. let $x_1b_{21} \notin E(G) \Rightarrow$ By **Claim 1**, $x_1c_{11} \notin E(G) \Rightarrow x_1b_{1i} \in E(G)\ \forall\ i$. If $\exists\ x_2, x_3$ s.t. $x_2c_{11}, x_3c_{11} \in E(G)$, then $x_ib_{1j} \notin E(G)$ for $\forall\ j$ and $i = 2, 3 \Rightarrow \forall\ t > 1$, $x_mb_{2t} \notin E(G)$ for some $m \in \{2, 3\}$ (else $\langle b_{2t}, x_2, y_1, x_3, b_{11} \rangle = 2K_1+(K_2 \cup K_1)$). W.l.g. let $x_2b_{22} \notin E(G)$. But then $\langle x_2, x_1, b_{22}, y_1, b_{21} \rangle = 2K_1+(K_2 \cup K_1)$). Hence w.l.g. let $x_2c_{11} \notin E(G) \Rightarrow$ By **A.1**, $x_2b_{21} \notin E(G)$. Then $|B_1| \leq 1$ (else $x_ib_{22} \in E(G)$ for $i = 1, 2$ and $\langle b_{22}, x_1, w, x_2, c_{11} \rangle = 2K_1+(K_2 \cup K_1)$). Again for $j > 1$, no two $b_{2i}$s have a common non-adjacent $x \in X$ (else if $xc_{11} \in E(G)$, then $\langle c_{11}, b_{2k}, v, b_{2t}, x \rangle = 2K_1+(K_2 \cup K_1)$ and if $xc_{11} \notin E(G)$, then $\langle b_{11}, b_{2k}, y_1, b_{2t}, x \rangle = 2K_1+(K_2 \cup K_1)$). Hence color pairs of corresponding non-adjacent vertices of $X$ and $B_2$-$b_{21}$; $(y_i, y_i'); (v, w)$ and $(b_{11}, c_{11})$ by same colors and remaining vertices of $X$ and $b_{21}$ by new colors to get $\chi \leq |X|+|Y|+1+2| \leq \omega+3$, a contradiction. Similarly it can be proved that $|C_2| \leq 1$.

**Case 2.1:** Every $b_{1i}$ is non-adjacent to a unique $c_{1i}$ and vice versa.
Then $|B| \leq 4$ (else w.l.g. $y_1b_i \in E(G)$ for $1 \leq i \leq 3$ and by **Claim 5** say $b_i = b_{1i}$ for $1 \leq i \leq 2$. Then $c_{11}y_1 \notin E(G)$ and $\langle c_{11}, b_{12}, y_1, b_3, w \rangle = 2K_1+(K_2 \cup K_1)$). Thus as $\Delta \geq 2\omega-1 \geq 11$, $|X \cup Y \cup Y'| \geq 7$ and $|X| \geq 1$. Also $|B_1| \geq 3$ (else color vertices of $X$ by $|X|$ new colors and pairs $(v, c_{21}); (w, b_{21}); (b_{1i}, c_{1i}); (y_i, y_i')$ by same colors to get $\chi \leq |X|+|Y|+2+|B_1| \leq \omega+3$). W.l.g. let $y_1b_{11}, y_1b_{12} \in E(G) \Rightarrow y_1'b_{13} \in E(G)$, $y_1'b \in E(G)\ \forall\ b \in B-\bigcup_{i=1}^{2}b_{1i}$ and $y_1c \in E(G)\ \forall\ c \in C-\bigcup_{i=1}^{2}c_{1i} \Rightarrow |C_2| = 0$ (else as $|B| \geq |C|$ and $|B_1| = |C_1|$, $|B_2| > 0$ and $\langle b_{21}, c_{13}, y_1, c_{21}, v \rangle = 2K_1+(K_2 \cup K_1)$). Now $|C_1| \neq 4$ (else $x_i$ is non-adjacent to at the most 4 vertices and deg $x_i >$ deg $v = \Delta$), $|B_2| > 0$ (else color vertices of $X$ by $|X|$ new colors and pairs $(v, w); (b_{1i}, c_{1i}); (y_i, y_i')$ by same colors to get $\chi \leq |X|+|Y|+1+|B_1| \leq \omega+3$) and $x_ib_{21} \notin E(G)$ (else deg $x_i > \Delta$). Then color vertices of $X-x_i$ by $|X|-1$ new colors and pairs $(x_i, b_{21}); (v, w); (b_{1i}, c_{1i}); (y_i, y_i')$ by same colors to get $\chi \leq |X|+|Y|+1+3 \leq \omega+3$, a contradiction.

**Case 2.2:** Some vertex $B_1$ or $C_1$ be non-adjacent to two vertices of $C_1$ or $B_1$.
Then by **Claim 4**, $|Y| = 1$ and by **Claim 2**, $b_{1i}y_j, c_{1m}y_n' \in E(G)$ for $\forall\ i, j, m, n$. Then $\forall\ i$, $b_{1i}c_{1j} \in E(G)$ for at the most one $j$ and vice versa. Also by **II**, $|X| > 0$. Let $X_1 = \{x_{1i} \in X / x_{1i}b_{1k} \notin E(G)$ where $b_{1k}c_{1m} \in E(G)$ for some $k, m\}$ and $X_2 = X-X_1$. Then $|X_1| > 0$ (else let $X' = \{x_{1i} \in X / x_{1i}b_{1t} \notin E(G)$ where $b_{1t}c_{1m} \notin E(G)\ \forall\ m\}$ and $X'' = X-X'$. Then $\Delta+s+1 = |X_1 \cup C_1 \cup y_1' \cup w|+|X'' \cup B_1 \cup y_1 \cup v|+|B_2 \cup C_2| \leq 2\omega+2) \Rightarrow |X_2| = 0$. Thus every $x_k$ is non-adjacent to some $b_{1i}$ where $b_{1i}c_{1j} \in E(G)$ for some $j \Rightarrow |B_1|, |C_1| \geq 2$ as $b_{1i}c_{1k}, b_{1t}c_{1j} \notin E(G)$ for some $k, t$. Now no two $x', x'' \in X$ are non-adjacent to the same $b_{1i}$ where $b_{1i}c_{1j} \in E(G)$ for some $j$. Let if possible $x'b_{1i}, x''b_{1i} \notin E(G)$. Then $xb_{1i} \notin E(G)\ \forall\ x \in X$ (else let $xb_{1i} \in E(G)$ and $xb_{1n} \notin E(G) \Rightarrow \langle x, x', b_{1n}, x'', b_{1i} \rangle = 2K_1+(K_2 \cup K_1)$). Also

$|B_2| = |C_2| = 1$ (else $\Delta+s+1 = |X\cup(C_1-c_{1j})\cup y_1'\cup w|+|B_1\cup y_1\cup v|+|c_{1j}\cup B_2\cup C_2| \le 2\omega+2$) and as $s > 1$ and $|Y| = 1$ $\exists$ x ∈ X s.t. $xb_{21}$ ($xc_{21}$) ∉ E(G) $\Rightarrow$ by **1.7** of **Lemma 1**, $|B_1|, |C_1| \le 2$. Then color vertices of X by $|X|$ new colors and pairs $(y_1, y_1')$; $(v, w)$; and those of corresponding non-adjacent vertices of $B_1$ and $C_1$ by same colors to get $\chi \le |X|+|Y|+1+|B_1| \le \omega+3$, a contradiction. Hence if $xb_{1i}$, $x'b_{1j}$ ∉ E(G), then i≠j.

W.l.g. let $|B_1| \le |C_1|$ and $(b_{1i}, c_{1i})$ be s.t. $b_{1i}c_{1i}$ ∈ E(G) for i = 1,.., q. Then $|X| \le q \le |B_1|, |C_1| \le \omega$-2. Now color the pairs $(x_{2r+1}, b_{2r+1}), (x_{2r+2}, c_{2r+2})$ ($1\le 2r+1\le|X|$ if $|X|$ is odd and $1\le 2r+2\le|X|$ otherwise), $(b_t, c_{t+1})$ ($|X|<t\le$mod $|B_1|$), $(v, w), (y, c_2), (z, b_2)$ by same colors and remaining vertices of $|C_1|$ by different new colors to get $\chi \le |X|+|C_1|-\dfrac{(|X|-1)}{2}+2 \le \omega+\dfrac{|X|+1}{2} \le \dfrac{3\omega}{2}$, a contradiction.

**Case 3:** $\forall$ **partitioning pair** $(v, w)$ with deg $v = \Delta$, $<B\cup C>$ is complete.

**Case 3.1:** $|C| > 1$.
If $\exists$ c, c' s.t. zc, zc' ∈ E(G), then by **1.7** of **Lemma 1**, $<y_1'\cup B>$ is complete. Also as $|B| \ge |C|, |B| \ge 2$ and $y_1'\cup C>$ is complete. But then deg $y_1' > \Delta$, a contradiction. Hence w.l.g. let $y_1'c, y_1c'$ ∈ E(G). Clearly $|C| = 2$ and also $|B| = 2$ with $y_1'b, y_1b'$ ∈ E(G). Then color pairs $(c, y_1); (b', y_1'); (v, c'); (w, b); (y_i, y_i')$ ($2\le i\le|Y|$) by same colors and vertices of X by $|X|$ new colors to get $\chi \le |X|+|Y|+3 \le \omega+2$, a contradiction.

**Case 3.2:** $|C| = 1$.
W.l.g. let $<Y'\cup c>$ be complete. By **1.5** of **Lemma 1**, every b ∈ B has at the most one vertex from Y adjacent. Let B' = {b ∈ B/ $<b\cup Y'>$ is not complete} and B" = B-B'.
We have
(a) As $s \ge 2$, every $y_i'$ is non-adjacent to some vertex of B' but by **1.7** of **Lemma 1**, $y_i'$ is not non-adjacent to more than one vertex of B'. Thus every vertex of B' is non-adjacent to one $y_i'$ and vice versa $\Rightarrow |B'| = |Y'|$. Also deg $y_i' = \Delta$ $\forall$ i.
(b) xc ∉ E(G) $\forall$ x ∈ X (else as $s \ge 2$, $xb_i$ ∉ E(G) for say i = 1, 2 and $<c, b_1, v, b_2, x> = 2K_1+(K_2\cup K_1)$)
(c) By **Claim** 1, xb' ∉ E(G) $\forall$ b' ∈ B'
(d) By **1.7** of **Lemma** 1, every vertex of X is adjacent to at the most one vertex of B"
(e) Every vertex of B" is adjacent to at the most one vertex of X (else $<b", x_1, w, x_2, c> = 2K_1+(K_2\cup K_1)$).
(f) If t = min{|B"|+1, |X|}, then $\exists$ $(x_i, u_i)$ where $x_i$ ∈ X, $u_i$ ∈ {B"∪c} s.t. $x_iu_i$ ∉ E(G). This can be proved by induction on t using (d) and (e).
(g) If $y_1'b_1'$ ∉ E(G), then roles of X and B" are interchanged when instead of $(v, w)$ $(y_1', b_1')$ is considered as a **partitioning pair**.
(h) $<X\cup Y\cup v>, <B"\cup Y'\cup c>, <Y'\cup c\cup w>$ are all complete.

By (g) w.l.g. let $|X| \ge |B"|$. If $|X| > |B"|$ and by (f) $(c, x_1), \bigcup_{i=2}^{|B"|}(x_i, b_i")$ be s.t. $x_ib_i"$ ∉ E(G), then color pairs $(x_1, c); (x_i, b_i")$ ($2\le i\le|B"|$); $(x_j, b_{|Y|-|X|+j}")$ ($|B"|+2\le j\le|X|$); $(b_k', y_k'), (b_{k+1}', y_k)$, ($1\le k\le \lfloor\dfrac{|Y|-|X|+|B"|+1}{2}\rfloor$) (if $|Y|-|X|+|B"|+1$ is odd then $y_k$ has no pair for k=$\lceil\dfrac{|Y|-|X|+|B"|+1}{2}\rceil$); $(y_m', y_m)$ ($\lceil\dfrac{|Y|-|X|+|B"|+1}{2}\rceil+1\le k\le|Y|$); $(v, w)$ by same colors to get $\chi \le |X|+2\lceil\dfrac{|Y|-|X|+|B"|+1}{2}\rceil+|Y|-\lceil\dfrac{|Y|-|X|+|B"|+1}{2}\rceil+1\le \dfrac{3|Y|+|X|+|B"|+4}{2} = \dfrac{(|Y|+|X|+1)+(|Y|+|B"|+1)+|Y|+2}{2} \le \dfrac{3\omega}{2}$ by (h), a contradiction. Next let $|X|=|B"|$. If $|X| > 0$, then as before we arrive at a contradiction. Finally if $|X| = 0$, then color pairs $(b_k', y_k'), (b_{k+1}', y_k)$,

$(1\leq k\leq \lceil \frac{|Y|}{2} \rceil)$; $(y_m', y_m)$ $(\lceil \frac{|Y|}{2} \rceil +1\leq k\leq |Y|)$; (v, w) by same colors. If |Y| is odd then color pair $(c, y_k)$ by the same color for $k=\lceil \frac{|Y|}{2} \rceil$) and otherwise color c by a new color to get $\chi \leq \frac{3|Y|+4}{2} < \frac{3\omega}{2}$ by (h), a contradiction.

This proves the **Main Result**.

**Examples of** $\{3K_1, 2K_1+(K_2\cup K_1)\}$-free graphs with tight $\chi(G)$.

1. If $\omega=2r$, then let $G \sim \sum_{i=1}^{r} C_5$, $\chi=3r$.

2. If $\omega=2r+1$, then let $G \sim \sum_{i=1}^{r} C_5 + W_6$, $\chi=3r+1$

3. If $\omega=5$, then let G be as follows: $V(G) = \{v, w, \bigcup_{i=1}^{3} y_i, \bigcup_{i=1}^{3} y_{i'}, \bigcup_{i=1}^{4} b_i, \bigcup_{i=1}^{4} c_i \}$.

   E(G) is defined as below:

| Vertex | Non-adjacent to vertices | Vertex | Non-adjacent to vertices |
|---|---|---|---|
| v | w; $c_i$ $1\leq i\leq 4$ | w | v, $b_i$ $1\leq i\leq 4$ |
| $b_1$ | w, $c_1$, $y_1'$, $y_2'$, $y_3'$ | $c_1$ | v, $b_1$, $y_1$, $y_2$, $y_3$ |
| $b_2$ | w, $c_2$, $y_1$, $y_2$, $y_3'$ | $c_2$ | v, $b_2$, $y_1'$, $y_2'$, $y_3$ |
| $b_3$ | w, $c_3$, $y_1$, $y_2'$, $y_3$ | $c_3$ | v, $b_3$, $y_1'$, $y_2$, $y_3'$ |
| $b_4$ | w, $c_4$, $y_1'$, $y_2$, $y_3$ | $c_4$ | v, $b_4$, $y_1$, $y_2'$, $y_3'$ |
| $y_1$ | $y_1'$, $c_1$, $b_2$, $b_3$, $c_4$ | $y_1'$ | $y_1$, $b_1$, $c_2$, $c_3$, $b_4$ |
| $y_2$ | $y_2'$, $c_1$, $b_2$, $c_3$, $b_4$ | $y_2'$ | $y_2$, $b_1$, $c_2$, $b_3$, $c_4$ |
| $y_3$ | $y_3'$, $c_1$, $c_2$, $b_3$, $b_4$ | $y_3'$ | $y_3$, $b_1$, $b_2$, $c_3$, $c_4$ |

**Edge Set E(G)**

Then G is $\{3K_1, 2K_1+(K_2\cup K_1)\}$-free, 10-regular and $\chi(G) = 8$.